\documentclass[english,reqno]{amsart}

\usepackage{babel}
\usepackage{amstext}
\usepackage{amsmath}
\usepackage{amsfonts}
\usepackage{latexsym}
\usepackage{ifthen}
\usepackage{xypic}
\usepackage{paralist}

\makeatletter
\@addtoreset{equation}{section}

\makeatother

\theoremstyle{plain} 
\newtheorem{theorem}{\noindent\bf Theorem}[section]
\newtheorem{lemma}[theorem]{\noindent\bf Lemma}
\newtheorem{corollary}[theorem]{\noindent\bf Corollary}
\newtheorem{proposition}[theorem]{\noindent\bf Proposition}

\theoremstyle{definition} 
\newtheorem{remark}[theorem]{\noindent\bf Remark}

\newcommand{\bN}{{\mathbb N}}
\newcommand{\Aut}{\mathop{\rm Aut}}
\newcommand{\nsg}{\unlhd}
\newcommand{\gnsg}[1]{\langle\langle #1 \rangle\rangle}

\newcommand{\bT}{{\mathbb T}}
\newcommand{\bG}{\tilde{G}}
\newcommand{\hG}{\hat{G}}

\renewcommand{\mathbb}{\mathbf}
\newcommand{\zero}{^{\circ}}

\title[] {The fundamental group of quotients of a product of curves}
\author{Thomas Dedieu and Fabio Perroni}

\date{\today}

\begin{document}

\maketitle

\hspace{1,3cm}{\em This article is dedicated to the memory of Fritz Grunewald}

\begin{abstract}
We prove a structure theorem for the fundamental group of the quotient $X$
of a product of curves by the action of a finite group $G$, hence for that of
 any resolution of the singularities of $X$.
 \end{abstract}

\section{Introduction}

The study of varieties isogenous to a product of curves was initiated
by Catanese in \cite{fabrizio},  inspired by a construction of
Beauville. These varieties 
are quotients of a product of smooth projective curves
$C_1\times\cdots\times C_n$ by the free action of a finite group
$G$.

Much of the work in this area has been focused in the $n=2$
case. Surfaces 
isogenous to a product of curves provide a wide class of surfaces
quite manageable to work with, since they are determined by discrete
combinatorial data. They were used successfully to address various
questions (see e.g. the survey paper \cite{survey}), and in particular
to obtain substantial information about various moduli spaces of
surfaces of general type (see e.g. \cite{newsurfaces,BCG,BCGP}).

In the case of a variety isogenous to a product, the action of $G$ is
free, and $X:=(C_1\times\cdots\times C_n)/G$ is smooth. Furthermore,
we have the following natural description of the fundamental group of
$X$. 

\begin{proposition}
\label{prop:initial-result}
\cite{fabrizio}
If $X:=(C_1\times\cdots\times C_n)/G$ is the quotient of a product of
curves by the free action of a finite group, then the fundamental
group of $X$ sits in an exact sequence 
\begin{equation}
\label{natural-exseq}
1 \to \Pi_{g_1}\times\cdots\times\Pi_{g_n} \to \pi_1(X) \to G \to 1,
\end{equation}
where each $\Pi_{g_i}$ is the fundamental group of $C_i$. This
extension, in the unmixed case where  each
$\Pi_{g_i}$ is a normal subgroup,  is determined by the associated
maps $G\to \mathrm{Out}\left(\Pi_{g_i}\right)$ to the respective
Teichm\"uller modular groups. 
\end{proposition}

\medskip
In the recent paper \cite{BCGP}, Bauer, Catanese, Grunewald and
Pignatelli prove that a similar statement still holds under weaker
assumptions.

\begin{theorem}
\label{thmbcgp}
\cite[Thm. 0.10 and Thm. 4.1]{BCGP}
Assume that $G$ acts faithfully on each curve
$C_i$ as a group of automorphisms, and let $X:=(C_1\times \cdots
\times C_n)/G$ be the (possibly singular) quotient by the diagonal
action of $G$. Then the fundamental group $\pi_1(X)$ has a normal
subgroup of finite index isomorphic to the product of $n$ surface
groups. We call $G'$ the quotient group. 
\end{theorem}

Here, by a surface group we mean a group isomorphic to the fundamental
group of a compact Riemann surface.
Note that, unlike in Proposition \ref{prop:initial-result}, the surface
groups in Theorem \ref{thmbcgp} above are not 
necessarily isomorphic to the fundamental groups of the curves
$C_1,\ldots,C_n$, and furthermore that the corresponding quotient $G'$ of
$\pi_1(X)$ is not necessarily isomorphic to $G$.

The first step of the proof of Theorem \ref{thmbcgp} consists in
showing that $\pi_1(X)$ is isomorphic to the quotient of the fibre
product $\bT:=\bT_1\times_G\cdots\times_G\bT_n$ of $n$ orbifold surface
groups (see Subsection \ref{notations}) by its torsion
subgroup $\mathrm{Tors}(\bT)$. Whereas this first part rests upon
geometrical considerations,
the rest of the proof relies on  an abstract group theoretic argument
showing that this quotient necessarily contains a normal subgroup as
described in Theorem \ref{thmbcgp}.
In particular, the relation occurring between the groups $G$ and $G'$
is not well understood. 

Using a suitable resolution of the singularities of
$X$, Bauer, Catanese, Grunewald and Pignatelli show in addition that
the fundamental group of any resolution  $Y$ of  $X$ is
isomorphic to the fundamental group of $X$, so  that the same description holds for $\pi_1(Y)$.

In \cite{BCGP}, as an important application of Theorem \ref{thmbcgp},
many new families of algebraic surfaces $S$ of general type  
with $p_g(S)=0$ are constructed, and several new examples of groups
are realized as the fundamental group of an algebraic surface $S$ of
general type with $p_g(S)=0$. This increases notably our
knowledge on algebraic surfaces. In fact the authors consider and
classify all the 
 surfaces whose canonical models arise as quotients
$X:=\left(C_1\times C_2\right)/G$ 
of the product of two curves of genera $g(C_1),g(C_2)\geq 2$ by
the action of a finite group $G$ such that $p_g(X)=q(X)=0$.

\medskip

In the present paper, we drop the assumption that the actions of
$G$ on $C_1,\ldots,C_n$ are faithful. We obtain the following expected
strengthening of Theorem \ref{thmbcgp}.

\begin{theorem}
\label{theo:main-intro}
Let $C_1,\ldots,C_n$ be smooth projective curves, and let $G$ be a
finite group acting on each $C_i$ as a group of automorphisms.
Then the fundamental group of the quotient
$X:=C_1\times\cdots\times C_n/G$ by the diagonal action of $G$ has a
normal subgroup of finite index that is isomorphic to the product of
$n$ surface groups.
\end{theorem}

This result should allow in the future the realization of interesting
groups as fundamental groups of higher dimensional  algebraic
varieties, following the method developed in \cite{BCGP} for
surfaces. Notice that, in the case where the $G$-actions are faithful,  
$X$ can only have isolated cyclic-quotient singularities, while
if the actions are not faithful, then the singular locus of $X$ can
have components of positive dimension, and the singularities 
are abelian-quotient singularities.

Again, one shows that any desingularization of the quotient $X$ has a
fundamental group isomorphic to that of $X$. This time however, we have
 to rely on a strong result of Koll\'ar \cite{kollar}.

The proof follows then closely the one of Theorem \ref{thmbcgp} of
\cite{BCGP}.
The main new difficulty one has to overcome is to find a
natural counterpart to the fibered product $\bT_1\times_G\cdots\times_G
\bT_n$, acting discontinuously on the product
$\tilde{C}_1\times\cdots\times \tilde{C}_n$ of the universal covers of
$C_1,\ldots,C_n$. After that similar group theoretic arguments work with
some slight modifications.

It has already been observed in \cite{fabrizio} that Theorem
\ref{theo:main-intro} follows directly from Theorem \ref{thmbcgp} when
$n=2$, by performing the quotient $\left(C_1\times C_2\right)/G$ in
successive steps. For $n>2$ however, this procedure
does not apply. 

The paper is organized as follows. In Section \ref{sec:background}, we
fix notations and collect some basic facts about group actions on
compact Riemann surfaces. 
Section \ref{sec:main} is devoted to the proof of Theorem
\ref{theo:main-intro}: the proof itself is given in Subsection
\ref{subsec:proofmain}, using intermediate results proven in
Subsections \ref{subsec:gt} and \ref{subsec:extension}.

\medskip
{\it Acknowledgements.} We are indebted to Prof. Ingrid Bauer
for suggesting this problem to us and for many useful
conversations. 
We express our thanks to  Profs. Fabrizio Catanese and Keiji
Oguiso  for their interest and warm encouragement. In particular, we
thank Prof. Oguiso for bringing the paper \cite{kollar} to our
attention. 
We also wish to thank Wenfei Liu and Matteo Penegini for
worthy conversations on the subject.

This work started as the first-named author held a post-doctoral
position at the Universit\"at Bayreuth, supported by the
DFG-Forschergruppe 790 ``Classification of Algebraic Surfaces and
Compact Complex Manifolds''. He is grateful to this institution for
providing excellent working conditions.

\section{Notations and basic results}
\label{sec:background}

\subsection{Notations}\label{notations}

We work over the field of complex numbers $\mathbb{C}$. 

Let $G$ be a group, and consider a subset $H \subset G$. We write
$H\leq G$ when $H$ is a subgroup of $G$, and 
$H\unlhd G$ when $H$ is a normal subgroup of $G$. If $A \subset G$ is
any subset, then $\gnsg{A}_G$ denotes the normal subgroup of $G$
generated by $A$.

\medskip
Let $g$ be a non negative integer. 
We call $\Pi_g$ the \textit{surface group} of genus $g$, defined as 
\begin{align*}
\Pi_g := \left\langle \vphantom{\prod\nolimits_{i=1}^g} \, a_1,
  b_1,\dots , a_g, b_g\, \right| \,  
\left. \prod\nolimits_{i=1}^g [a_i,b_i] = 1 \,
\right\rangle\, .
\end{align*} 
This is the fundamental group of compact Riemann surfaces of genus
$g$. 
On the other hand, letting in addition $m_1,\ldots,m_r$ be positive
integers, we denote by $\mathbb{T}(g; m_1,\dots , m_r)$ the
\textit{orbifold surface group} of signature $(g;
m_1,\dots , m_r)$, defined as
\begin{multline*}
\mathbb{T}(g; m_1,\dots, m_r):= \left\langle
  \vphantom{\prod\nolimits_i^g} \, a_1,
  b_1,\dots, a_g,b_g, c_1, \dots , c_r \, \right| \, \\ 
\left. c_1^{m_1}= \cdots = c_r^{m_r} = 
\prod\nolimits_{i=1}^g [a_i ,b_i]\cdot c_1\cdot \, \cdots \,
\cdot c_r=1\, \right\rangle. 
\end{multline*}
It is obtained from the fundamental group of the complement of a
set of $r$ distinct points in a compact Riemann surface of genus $g$,
by quotienting by the normal subgroup generated by
$\gamma_1^{m_1},\ldots,\gamma_r^{m_r}$, where each $\gamma_i$ is a
simple geometric counterclockwise loop around the $i$-th removed
point.

Let $G$ be a finite group. An \textit{appropriate orbifold
  homomorphism} 
is a surjective homomorphism 
$\varphi : \mathbb{T}(g; m_1,\dots , m_r)\rightarrow G$
such that $\varphi (c_i)$ has order $m_i$
for $i=1,\ldots,r$.

\medskip
The action of a group $G$ as a group of homeomorphisms on a
topological space $X$ is said to be \textit{discontinuous} if the
following two conditions are satisfied: 
\begin{inparaenum}[(i)]
\item the stabilizer of each point of $X$ is finite;
\item each point $x\in X$ has a neighbourhood $U$ such that $g(U)\cap
  U=\emptyset$, for each $g\in G$ such that $gx\not=x$. 
\end{inparaenum}

\subsection{Basic results}

The following result is essentially a reformulation of Riemann's
existence theorem (see \cite[Thm. 2.1]{BCGP}).

\begin{theorem}
A finite group $G$ acts faithfully as a group of automorphisms on a
compact Riemann surface of genus $g$ if and only if there are natural
numbers $g', m_1,\ldots , m_r$ and an appropriate orbifold
homomorphism  
\[
\varphi : \mathbb{T}(g';m_1,\ldots , m_r)\rightarrow G
\]
such that the Riemann-Hurwitz relation holds:
$$
2g - 2 = |G|\left( 2g' - 2 + \sum_{i=1}^r \left( 1-
    \frac{1}{m_i}\right)\right). 
$$
\end{theorem}

\begin{remark}
\label{rem1}
As already remarked in \cite{BCGP}, under the above hypotheses, $g'$
is the geometric genus of $C':=C/G$, and $m_1,\dots , m_r$ are the
branching indices at the branching points of the $G$-cover $p \colon
C\rightarrow C'$. 
The appropriate orbifold homomorphism  $\varphi$ is induced by
the monodromy of the Galois \'etale $G$-covering $p\zero \colon
C\zero \to {C'}\zero$ induced by $p$, where ${C'}\zero$ is the Riemann 
surface obtained from $C'$ by removing the branch points of $p$, and
${C}\zero := p^{-1}({C'}\zero)$.
In particular, $\varphi(c_i)$ generates the stabilizer of the
corresponding ramification point. 

Furthermore, the kernel of   $\varphi$ is isomorphic to the
fundamental group $\pi_1(C)$, and  the action of $\pi_1(C)$ on the
universal cover $ \tilde{C}$ of $C$ extends to a discontinuous  action
of $\mathbb{T} := \mathbb{T}(g';m_1,\dots , m_r)$. 
Let $u\colon \tilde{C}\to C$ be the covering map. It is
$\varphi$-equivariant, and $C/G \,
\cong \, \tilde{C}/\mathbb{T}\, .$ 

\medskip

We now give two elementary facts that will be used in the following.

\begin{lemma}\label{l12}
\begin{inparaenum}[\textnormal{(\roman{enumi})}]
\item Let $x\in \tilde{C}$. Then the restriction of $\varphi$ to the
  stabilizer ${\rm St}_x$ of $x$ (with respect to the action of
  $\mathbb{T}$ on $\tilde{C}$) is injective. \\
\item Let $t\in {\rm St}_x$. Then $t$ is
  conjugated to $c_i^m$, for some $i\in \{1,\ldots, r\}$ and $m\in
  \bN$. 
\end{inparaenum}
\end{lemma}

\begin{proof} The $\pi_1(C)$-action on $\tilde{C}$ is free, so
  $\pi_1(C)\cap {\rm St}_x =\{1\}$. This yields (i), because
  $\pi_1(C)$ is the kernel of $\varphi$.

To prove (ii), let $y=u(x)$. If $t=1$, then the result is clear. Else,
there exists an integer $i\in\{1,\ldots,m\}$ and a point $x'\in
u^{-1}(y)$ that is fixed by $c_i$. It then follows from (i) that 
${\rm St}_{x'}=\langle c_i \rangle$, hence that $t$ is conjugated to a
power of $c_i$. 

\end{proof}

\end{remark}
 
\section{Main Theorem}
\label{sec:main}

\noindent
The main result of the paper is the following
\begin{theorem}\label{mainthm}
Let $C_1,\ldots, C_n$ be compact Riemann surfaces, and
let $G$ be a finite group that acts as a group of
automorphisms on each $C_i$.
We consider the quotient of the product $C_1\times \cdots \times C_n$
by the diagonal action of $G$.
Then there is a normal subgroup of finite index $\Pi$ in the
fundamental group 
$$
\pi_1 \left( \frac{C_1\times \cdots \times C_n}{G}\right),
$$ 
such that $\Pi$ is isomorphic to the product of $n$ surface groups.  
\end{theorem} 

The proof of this theorem follows closely  that of
\cite[Thm. 4.1]{BCGP}, and is given in the next subsections.
Before we move on to this proof, let us give the following important
consequence of Theorem \ref{mainthm}.

\begin{corollary}\label{thm2}
Let $C_1,\ldots, C_n$ and $G$ be as in the statement of
Theorem \ref{mainthm}, and let $Y$ be a resolution of the
singularities of $X:=(C_1\times \cdots \times C_n)/G$. 
Then, the fundamental group  of $Y$ has a  normal subgroup of finite
index  isomorphic to the product of $n$ surface groups.
\end{corollary}

\begin{proof} By Theorem \ref{mainthm}, it is enough to show that the
natural morphism 
\[
f_\ast : \pi_1(Y) \longrightarrow \pi_1(X)
\]
induced by the resolution  $f: Y \rightarrow X$ is an isomorphism.

This follows directly from \cite[Sec. 7]{kollar}: since $X$ is normal
and only has quotient singularities, $Y$ is locally simply connected
by \cite[Thm. 7.2]{kollar}, hence $f_*$ is an isomorphism by
\cite[Lem. 7.2]{kollar}. 

\end{proof}

\subsection{Proof of the main theorem}
\label{subsec:proofmain}

For $i=1,\ldots,n$, we let
\[
{K}_i ={\ker}\left( G\to \Aut(C_i) \right)\quad
\mbox{and}\quad H_i=G/K_i \, , 
\]
where $ G\to {\rm Aut}(C_i)$ is the morphism associated to the action
of $G$ on $C_i$. We call $p_i$ the projection $G \to H_i$.
Now $H_i$ acts faithfully on $C_i$, so we have (see Remark \ref{rem1})
a short exact sequence
\begin{equation}
\label{ret}
1\to \pi_1(C_i) \to \bT_i \xrightarrow{\varphi_i}
H_i \to 1,
\end{equation}
where $\bT_i$ is an orbifold surface group, and $\varphi_i$ is an
appropriate orbifold homomorphism. Let $\Sigma_i:=G\times_{H_i} \bT_i$
be the fibered product corresponding to the Cartesian diagram
\[ \xymatrix{
\Sigma_i \ar[d]_{\psi_i} \ar[r] \ar@{}[dr]|{\Box} & \bT_i
\ar[d]^{\varphi_i} \\ 
G \ar[r]_{p_i} & H_i.
} \]
We call $\psi_i:\Sigma_i \to G$ the projection on the first
factor. Pulling-back \eqref{ret} by $p_i:G\to H_i$, we obtain a
short exact sequence 
\begin{equation}
\label{ret2}
1\to \pi_1(C_i) \to \Sigma_i \xrightarrow{\psi_i}
G \to 1,
\end{equation}
where the left-hand side map is $\gamma \in \pi_1(C_i) \mapsto
(1,\gamma) \in \Sigma_i$. 

Next, we define $\bG:=\Sigma_1 \times_G \cdots \times_G \Sigma_n$ as
the fibered product correponding to the Cartesian diagram below.
\[ \xymatrix@R=15pt@C=10pt{
&& \bG \ar[1,-2] \ar[1,-1] \ar[1,1] \ar[1,2] && \\
\Sigma_1 \ar[1,2] & \Sigma_2 \ar[1,1] & \ldots & \Sigma_{n-1}
\ar[1,-1] & \Sigma_n \ar[1,-2] \\
&& G &&
} \]
Let $\Delta: G \to G \times \cdots \times G$ be the diagonal
morphism. Then $\bG$ can also be seen as the fibered product
$G\times_{(G\times\cdots\times
  G)} (\Sigma_1\times \cdots \times \Sigma_n)\to G$ with respect to
the two morphisms $\Delta$ and $(\psi_1,\ldots,\psi_n)$.
Therefore, the pull-back by $\Delta$ of the product of the $n$ exact 
sequences \eqref{ret2} for $i=1,\ldots,n$ is a short exact
sequence
\begin{equation}
\label{ret-final}
1 \rightarrow \prod\nolimits_{i=1}^n \pi_1(C_i)\to  \bG
\xrightarrow{\Psi} G \to 1,
\end{equation}
where $\Psi$ is the first projection $G\times_{(G\times\cdots\times
  G)} (\Sigma_1\times \cdots \times \Sigma_n)\to G$.

Now we have the following, coming from the fact that $\bG$ acts
discontinuously on the universal cover of $C_1\times \cdots \times
C_n$.

\begin{proposition}\label{prop1}
Let $\bG ' \unlhd \bG$ be the normal subgroup of $\bG$ generated by
those elements which have non-empty fixed-point set. 
Then
$$
\pi_1 \left( \frac{C_1\times\cdots\times C_n}{G}\right)
\quad \cong \quad \frac{\bG}{\bG '}. 
$$
\end{proposition}

\begin{proof}
For $i=1,\ldots,n$, the action of $\bT_i$ on the universal covering
$\tilde{C_i}$ of $C_i$ 
(see Remark \ref{rem1}) induces an action of $\Sigma_i$ on
$\tilde{C_i}$ via the projection of $\Sigma_i$ on its second factor
$\bT_i$. 
We obtain in this way an action of $\bG$ on the product $\tilde{C}_1
\times \cdots \times \tilde{C}_n$. 

This action is discontinuous: let ${\rm St}_x$ be the stabilizer of a
point $x \in \tilde{C}_1 \times \cdots \times \tilde{C}_n$ with
respect to the action of $\bG$. Then the same argument as that in the
proof of Lemma \ref{l12} shows  that $\Psi_{|{\rm St}_x}$ is
injective, from which it follows that ${\rm St}_x$ is finite because
$G$ is finite. On the other hand, 
condition (ii) in the definition of a discontinuous action is a
consequence of the fact that the $\bT_i$-actions are themselves
discontinuous. 

Then, the main theorem in \cite{armstrong} applies to our situation,
and gives a group isomorphism
\[
\pi_1\left( \frac{\tilde{C}_1 \times \cdots \times \tilde{C}_n}{\bG}
\right) \cong \frac{\bG}{\bG'}.
\]
Eventually, since the universal covering
$\mathcal{U}\colon 
\tilde{C}_1\times \cdots \times \tilde{C}_n \to C_1\times \cdots
\times C_n$ is $\Psi$-equivariant, we have an isomorphism
\[
\frac{C_1\times \cdots \times C_n}{G} \cong 
\frac{\tilde{C}_1\times\cdots\times\tilde{C}_n}{\bG},
\]
and the proposition follows. 

\end{proof}

\begin{remark}\label{tor}
The elements of $\bG$ which have fixed-points are precisely those
elements of finite order. Therefore $\bG ' $ is the torsion subgroup
of $\bG$.
\end{remark}

Now the proof of Theorem \ref{mainthm} relies on the following result,
the proof of which we postpone to Subsection \ref{subsec:extension}.

\begin{proposition}
\label{prop:extension}
The quotient $\bG/\bG'$ is an extension 
\[
1\to E \to \bG/\bG' \xrightarrow{\theta} T \to 1
\]
of a finite group $E$ by a group $T$ that is a finite-index
subgroup of a product of $n$ orbifold surface groups.
\end{proposition}

Using the results of \cite{GJ-ZZ08}, the latter fact enables one to
show that there is a finite index normal subgroup $\Gamma \nsg
\bG/\bG'$ that injects in $T$: 

\begin{lemma}
\label{prf}
Let $S$ be a group sitting in an exact sequence 
\begin{equation*}
1\to E \to S \to T \to 1,
\end{equation*}
where $E$ is a finite group, and $T$ is a finite index subgroup of a
product of $n$ orbifold surface groups. Then $S$ is residually
finite. In particular, there exists a finite index normal subgroup
$\Gamma \nsg S$ such that $\Gamma \cap E = \{1\}$.
\end{lemma}

\begin{proof} 
By \cite[Prop. 6.1]{GJ-ZZ08}, an extension of a finite
group by a group that is residually finite and good in the sense of
\cite{serre} is residually finite.
It therefore suffices to show that $T$ enjoys the two aforementioned
properties.

An orbifold surface group is residually finite. Therefore $T$ is
itself residually finite, being a finite index subgroup of a product
of orbifold surface groups.

By \cite[Lem. 3.2]{GJ-ZZ08}, it is enough to show that a product of
orbifold surface groups is good to prove that $T$ is good. But
\cite[Prop. 3.7]{GJ-ZZ08} tells us that an orbifold surface group is
good, and \cite[Prop. 3.4]{GJ-ZZ08} that a product of good groups is
good. 

\end{proof}

We are now in a position to complete the proof of our main theorem:

\begin{proof}[Proof of Theorem \ref{mainthm}]
Let $\bT_1\times\cdots\times\bT_n$ be a product of $n$ orbifold
surface groups containing $T$ as a finite index subgroup,
and let us consider $\Gamma\unlhd \bG/\bG'$ a normal subgroup of
finite index such that $E\cap\Gamma=\{1\}$. 
Then $\theta(\Gamma) \leq \bT_1\times\cdots\times\bT_n$ has finite
index. 

Now every orbifold surface group contains a surface group as a finite
index subgroup (see e.g. \cite{beardon}), so let $\Pi_i$ be a finite
index surface group in $\bT_i$ for each $i=1,\ldots,n$.

For each $i$, we consider
\[
\theta(\Gamma)_i:= \theta(\Gamma) \cap \left( \{1\} \times \cdots
\times \bT_i \times \cdots \times \{1\} \right)
\]
as a subgroup $\theta(\Gamma)_i \leq \bT_i$, and set
\[
\Pi_i':=\bigcap_{g\in\bT_i} g\left(\theta(\Gamma)_i\cap\Pi_i\right)
g^{-1},
\]
the biggest normal subgroup of $\bT_i$ contained in
$\theta(\Gamma)_i\cap\Pi_i$. Then $\Pi_i'$ has finite index in
$\Pi_i$, and thus is itself a surface group.
Eventually, $\Pi:=\Pi_1'\times\cdots\times \Pi_n'$ is a subgroup of
$\theta(\Gamma)$, which is normal and of finite index in
$T$. Therefore, $\theta^{-1}(\Pi)\cap
\Gamma$ is a normal subgroup of $\bG/\bG'$, with finite index, and
isomorphic to $\Pi$.

\end{proof}

\subsection{Results in group theory}
\label{subsec:gt}
In this subsection, we prove some technical results that are needed
for the proof of Proposition \ref{prop:extension}.

Let $\Sigma$ be any group, $R\unlhd \Sigma$ be a normal subgroup, and
$L\subset \Sigma$ be a subset. We define 
\begin{equation}
\label{def:N}
N(R,L):=\langle\langle \{ hkh^{-1}k^{-1}\, | \, h \in L\, , \, k\in R
\}\rangle \rangle_{\Sigma} 
\end{equation}
and
\begin{equation}
\label{def:hatsigma}
\hat{\Sigma}:=\hat{\Sigma}(R,L):=\Sigma/N(R,L).
\end{equation}
We call $\hat{R}$ and $\hat{L}$ the images of $R$ and $L$ respectively
by the projection $\Sigma \to \hat{\Sigma}$. There is an
isomorphism: $\hat{\Sigma}/\gnsg{\hat{L}}_{\hat{\Sigma}} \cong
\Sigma/\gnsg{L}_{\Sigma}$.
Notice also that $N(R,L)\unlhd R$ and $N(R,L)\unlhd
\gnsg{L}_{\Sigma}$, which implies that $\hat{R}$ is a normal subgroup
of $\hat{\Sigma}$.

\begin{lemma}
\label{finite}
If $R\unlhd \Sigma$ has finite index, and if $L\subset \Sigma$ is a
finite subset consisting of elements of finite order, then
$\gnsg{\hat{L}}_{\hat{\Sigma}}$ is finite.
\end{lemma}

\begin{proof} The subgroup $\langle\langle
  \hat{L}\rangle\rangle_{\hat{\Sigma}}$ 
is the image of $\langle\langle {L}\rangle\rangle_{{\Sigma}}$ under
the projection 
$\Sigma\to \hat{\Sigma}$. 
Since $R$ has finite index in $\Sigma$, and $L$ is finite,  it follows
that  $\langle\langle \hat{L}\rangle\rangle_{\hat{\Sigma}}$ is generated by
finitely many elements which are conjugated to those of
$\hat{L}$. Since the elements of $L$ have finite order, these
generators have finite order as well.

The center ${\rm Z}\left(\gnsg{\hat{L}}_{\hat{\Sigma}}\right)$ of
$\langle\langle \hat{L}\rangle\rangle_{\hat{\Sigma}}$ contains
$\hat{R}\cap \langle\langle 
\hat{L}\rangle\rangle_{\hat{\Sigma}}$, and hence has finite index in
$\gnsg{\hat{L}}$.  Now, by \cite[Lem. 4.6]{BCGP}, if a group $S$ is
generated by finitely many elements of finite order, and if its centre
has finite index in $S$, then $S$ is finite. From this we conclude
that $\langle\langle\hat{L}\rangle\rangle_{\hat{\Sigma}}$ is finite.  

\end{proof}

We now consider the particular case when $\Sigma$ is a group
constructed as in Subsection \ref{subsec:proofmain}:
$\Sigma=G\times_H\bT$, where $G$ is a finite group, $H$ is a quotient
of $G$, and $\bT$ is any group coming with a surjective morphism
$\varphi:\bT \rightarrow H$.

\begin{lemma}\label{sigma}
The projection on the second factor $q : \Sigma \to \bT$ induces a
morphism
\begin{equation}\label{chi}
\bar{q}: \frac{\Sigma}{\langle\langle L\rangle\rangle_{\Sigma}}\to
\frac{\bT}{\langle\langle q(L) \rangle \rangle_{\bT}}
\end{equation}
in a natural way. It is surjective, and has finite kernel.
\end{lemma}

\begin{proof} We have $q\left(\gnsg{L}_{\Sigma}\right)=
\gnsg{q(L)}_{\bT}$. The map $\bar{q}$ is therefore induced by
the composition $\Sigma \xrightarrow{q} \bT \to
\bT/\gnsg{q(L)}_{\bT}$, which  is clearly surjective. To prove the
finiteness of its kernel, notice that for any $(g,t)\in
q^{-1}\left(\gnsg{q(L)}_{\bT}\right)$, there exists $h\in G$ with
$(h,t)\in \gnsg{L}_{\Sigma}$, hence $gh^{-1}\in K:= \ker (G\to H)$.
It follows that $q^{-1}(\langle\langle q(L)\rangle\rangle_{\bT})=
K\gnsg{L}_{\Sigma}$, where $K$ is seen as a subgroup in $\Sigma$ via
the injection $k\in K \mapsto (k,1)\in\Sigma$. Eventually,
${\ker}(\bar{q})\cong {K} \subset G$, which is finite. 

\end{proof}

\subsection{Realization of the fundamental group as a suitable extension}
\label{subsec:extension}

In this subsection, we give a full proof of Proposition
\ref{prop:extension}. We use the basic results in group theory
established in Subsection \ref{subsec:gt} above.

For $i=1,\ldots,n$, we fix the following presentation for the orbifold
groups $\bT_i$ in \eqref{ret}:
\begin{multline*}
\bT_i= \left\langle
  \vphantom{\prod\nolimits_i^g} a_{i1}, b_{i1},\ldots, a_{ig_i},
  b_{ig_i}, c_{i1}, \ldots , c_{ir_i} \, \right| \, \\  
\left. c_{i1}^{m_{i1}} = \cdots  =  c_{ir_i}^{m_{ir_i}} = 
\prod\nolimits_{j=1}^{g_i} [a_{ij} , b_{ij}]\cdot c_{i1}\cdot \,
\cdots \, \cdot c_{ir_i} =1 \right\rangle,
\end{multline*}
and set $R_i = \pi_1(C_i)$.
We write the elements of $\bG$ as $(g,z_1,\ldots,z_n)$,
with $(g,z_i)\in \Sigma_i$ for $i=1,\ldots,n$.
Then we have:

\begin{lemma}
\label{ni}
For each $i=1,\ldots,n$, there exists a finite subset $\mathcal{N}_i
\subset \bG$, such that
\[ \gnsg{\mathcal{N}_i}_{\bG} = \bG', \]
and whose elements are of the form 
\[
(g,z_1d_1^{\ell_1}z_1^{-1},\ldots, d_i^{\ell_i},\ldots,
z_nd_n^{\ell_n}z_n^{-1})
\]
for some $g\in G$, some $d_j \in \{c_{j1},\ldots,c_{jr_j}\}$ and
$\ell_j \in \bN$ for $j=1,\ldots,n$, and some $z_j \in \bT_j$ for $j
\neq i$. 
\end{lemma}

\begin{remark}
As a direct consequence of Lemma \ref{ni}, if
\[(g,z_1d_1^{\ell_1}z_1^{-1},\ldots, d_i^{\ell_i},\ldots,
z_nd_n^{\ell_n}z_n^{-1}) \in \mathcal{N}_i\] for some $i$, then for any
$j \neq i$, there exists \[(h,y_1,\ldots,\widehat{y_j},\ldots,y_n) \in
\Sigma_1 \times_G \cdots \times_G \widehat{\Sigma_j} \times_G \cdots
\Sigma_n\]
(where a hat means that the corresponding factor is omitted), such
that 
\[
(hgh^{-1},y_1z_1d_1^{\ell_1}z_1^{-1}y_1^{-1},\ldots,
y_id_i^{\ell_i}y_i^{-1},\ldots, d_j^{\ell_j},\ldots,
y_nz_nd_n^{\ell_n}z_n^{-1}y_n^{-1}) \in \mathcal{N}_j.
\]
\end{remark}

\begin{proof}[Proof of Lemma \ref{ni}]
Let $s \in \bG$ be an element with non empty fixed-point set, and let
us fix $i \in \{1,\ldots,n\}$. 
By Lemma \ref{l12} (ii), $s$ writes
\[
s= (g,z_1d_1^{\ell_1}z_1^{-1},\ldots,z_nd_n^{\ell_n}z_n^{-1}),
\]
with notations as in the statement of the Lemma.
Obviously, one can find $h\in G$, and
$\zeta_j\in\bT_j$ for each $j\neq i$, 
such that 
$(h,\zeta_1,\ldots,z_i^{-1},\ldots,\zeta_n) \in \bG$, and therefore
$s$ is conjugated in $\bG$ to an element of type
\begin{equation}\label{forms}
(g',y_1d_1^{\ell_1}y_1^{-1},\ldots,d_i^{\ell_i},\ldots,
y_nd_n^{\ell_n}y_n^{-1}).
\end{equation}

Now we claim that there exists finite sets $A_j \subset \bT_j$,
$j=1,\ldots,\widehat{i},\ldots,n$, such that each element of $\bG$ as
in \eqref{forms} is conjugated in $\bG$ to some
\[
(g'',x_1d_1^{\ell_1}x_1^{-1},\ldots,d_i^{\ell_i},\ldots,
x_nd_n^{\ell_n}x_n^{-1}),
\]
with $x_j \in A_j$ for each $j\neq i$.
Then it is clear that one can build $\mathcal{N}_i$ as required.

To prove our claim, first note that if
$(g,z_1d_1^{\ell_1}z_1^{-1},\ldots,d_i^{\ell_i},\ldots,
z_nd_n^{\ell_n}z_n^{-1}) \in \bG$, then an $(n+1)$-uple 
$(g,\zeta_1d_1^{\ell_1}\zeta_1^{-1},\ldots,d_i^{\ell_i},\ldots,
\zeta_nd_n^{\ell_n}\zeta_n^{-1})$ corresponds to an element of the
fibered product $\bG$ if and only if for each $j\neq i$,
$\varphi_j(z_j^{-1}\zeta_j)$ belongs to the centralizer
$C_{H_j}(\varphi_j(d_j^{\ell_j}))$ of $\varphi_j(d_j^{\ell_j})$ in
$H_j$. 

Second, note that if $k_j \in R_j$ for some $j\neq i$, then
$(1,1,\ldots,k_j,\ldots,1) \in \bG$, and therefore any element
$(g,\ldots,z_jd_j^{\ell_j}z_j^{-1},\ldots)\in\bG$ is conjugated to 
\[(g,\ldots,(k_jz_j)d_j^{\ell_j}(k_jz_j)^{-1},\ldots)\in\bG.\]
Then our claim follows from the fact that for each $j \neq i$, $R_j
\nsg \varphi_j^{-1}\left(C_{H_j}(\varphi_j(d_j^{\ell_j}))\right)$ has
finite index.

\end{proof}

From now on, we let $\mathcal{N}_1,\ldots,\mathcal{N}_n$ be as in
Lemma \ref{ni}.

\begin{lemma}\label{lemma4.10}
For $i=1,\ldots,n$, if 
\[(g,z_1d_1^{\ell_1}z_1^{-1},\ldots,
d_i^{\ell_i},\ldots, z_nd_n^{\ell_n}z_n^{-1}) \in \mathcal{N}_i,\]
then for all $k_i\in R_i$, we have
\[(1,1,\ldots,d_i^{\ell_i}k_id_i^{-\ell_i}k_i^{-1},\ldots,1) \in
\bG'. \]
\end{lemma}

\begin{proof} Let $k_i\in R_i$. Then
  $\tilde{k}_i:=(1,1,\ldots,k_i,\ldots,1) \in \bG$, and our result
  follows from the equality 
\begin{multline*}
(1,1,\ldots,d_i^{\ell_i}k_id_i^{-\ell_i}k_i^{-1},\ldots,1) =\\
(g,z_1d_1^{\ell_1}z_1^{-1},\ldots,d_i^{\ell_i},\ldots,
z_nd_n^{\ell_n}z_n^{-1}) \,
\tilde{k}_i \,
(g,z_1d_1^{\ell_1}z_1^{-1},\ldots,d_i^{\ell_i},\ldots,
z_nd_n^{\ell_n}z_n^{-1})^{-1}\,
\tilde{k}_i^{-1},
\end{multline*}
and the fact that $\gnsg{\mathcal{N}_i}_{\bG}=\bG'$.

\end{proof}

For $i=1,\ldots,n$, we let $L_i \subset \Sigma_i$ be the image of
$\mathcal{N}_i$ by the projection $\bG \to \Sigma_i$. The first
projection $\psi_i:\Sigma_i \to G$ then induces an epimorphism
$\hat{\psi}_i:\hat{\Sigma}_i:=\hat{\Sigma}_i(R_i,L_i)\to G$ (see
Subsection \ref{subsec:gt} for a definition of $\hat{\Sigma}_i$).

Eventually, we let $\hG$ be the fibered product
\begin{equation}
\label{hG}
\hat{\Sigma}_1\times_G\cdots \times_G \hat{\Sigma}_n \cong 
G \times_{(G\times\cdots\times G)} (\hat{\Sigma}_1\times\cdots\times
\hat{\Sigma}_n),
\end{equation}
and we define a map $\Phi : \hG \to  \bG/\bG'$ 
by the formula
\begin{equation}
\label{Phi}
\Phi\left([s_1],\ldots,[s_n]\right) =
\left[(s_1,\ldots,s_n)\right],
\end{equation}
where $s_i \in \Sigma_i$ for each $i$ (here we see $\hG$ as contained
in $\hat{\Sigma}_1\times\cdots\times \hat{\Sigma}_n$, using its
description by the left-hand side of \eqref{hG} rather than by its
right-hand side, and similarly we see $\bG$ as
contained in $\Sigma_1\times\cdots\times \Sigma_n$). 
It is a consequence of Lemma \ref{lemma4.10} that $\Phi$ is
well-defined by \eqref{Phi}.
Now we have:

\begin{lemma}
The morphism $\Phi$ is surjective, and has finite kernel.
\end{lemma}

\begin{proof} The surjectivity follows at once from \eqref{Phi}.
On the other hand, an element $\left([s_1],\ldots,[s_n]\right) \in
\hG$ lies in $\ker \Phi$ if and only if $(s_1,\ldots,s_n) \in
\bG'$. This implies for $i=1,\ldots,n$ that $s_i \in
\gnsg{L_i}_{\Sigma_i}$, because $\bG'=
\gnsg{\mathcal{N}_i}_{\bG}$. We thus see that $\ker \Phi$, seen as
contained in $\hat{\Sigma}_1\times\cdots\times \hat{\Sigma}_n$, is
contained in 
$ \gnsg{\hat{L}_1}_{\hat{\Sigma}_1} \times \cdots \times
\gnsg{\hat{L}_n}_{\hat{\Sigma}_n}$,
which is finite by Lemma \ref{finite}.

\end{proof}

Next, we define a morphism 
\begin{equation}
\label{Theta}
\Theta :\ 
\prod_{i=1}^n \hat{\Sigma}_i \ \to\ 
\prod_{i=1}^n \frac{\Sigma_i}{\gnsg{L_i}}
\ \xrightarrow{(\bar{q}_1,\ldots,\bar{q}_n)}\ 
 \prod_{i=1}^n \frac{\bT_i}{\gnsg{q_i(L_i)}}
\end{equation}
as in Subsection \ref{subsec:gt}: the left-hand side map in
\eqref{Theta} is the product of the projections 
\[\hat{\Sigma}_i \to \hat{\Sigma}_i/\gnsg{\hat{L}_i}_{\hat{\Sigma}_i}
\cong \Sigma_i/\gnsg{{L}_i}_{{\Sigma}_i},\]
and the $\bar{q}_i$'s are induced by the second projections $q_i:
\Sigma_i \to \bT_i$ as in Lemma \ref{sigma}.
We have
\begin{equation}\label{*}
\ker\Phi \subset \gnsg{\hat{L}_1}_{\hat{\Sigma}_1} \times\cdots\times
\gnsg{\hat{L}_1}_{\hat{\Sigma}_n} 
\subset \ker\Theta, 
\end{equation}
and $\ker\Theta$ is finite by both Lemmas \ref{finite} and
\ref{sigma}. 

Let us set 
\[ T :=\, \Theta (\hG). \]
Notice that $\hG$ has finite index in $\prod_{i=1}^n \hat{\Sigma}_i$
because $G$ is finite, and therefore that $T$ has finite index in
$\prod_{i=1}^n \bT_i/\gnsg{q_i(L_i)}$. 
We have a short exact sequence 
\begin{equation}\label{e2}
1\to E_1 \to \hG \xrightarrow{\left.\Theta\right|_{\hG}}
T \to 1. 
\end{equation}
Clearly, $\ker\Phi \subset E_1$. Therefore, setting $E:=E_1/\ker\Phi$,
we obtain the following commutative diagram
\begin{equation}
\label{final-diag}
\begin{xy}
  \xymatrix{
& 1\ar[d] & 1 \ar[d] &&\\
1 \ar[r] & \ker\, \Phi \ar[d] \ar[r]^{=}&\ker\, \Phi
\ar[r]\ar[d]&1\ar[d]&&\\ 
1\ar[r]&E_1\ar[r]\ar[d]&
\hG\ar[d]_{\Phi}\ar[r]^{\left.\Theta\right|_{\hG}}&
T\ar[d]_{=}\ar[r]&1\\ 
1\ar[r]&E\ar[r]\ar[d]&\bG/\bG'\ar[d]\ar[r]^{\theta}&
T\ar[d]\ar[r]&1\\  
&1&1&1 
}
\end{xy}
\end{equation}
where $\theta$ is the morphism induced by $\left.\Theta\right|_{\hG}$
which makes the diagram commutative.

We then claim that the lower row of the diagram \eqref{final-diag} is
the short exact sequence we are looking for:
exactness follows from an easy diagram chase; 
the finiteness of $E$ follows from that of $E_1$;
eventually, each $\bT_i/\gnsg{q_i(L_i)}$ is an orbifold surface
group, because $q_i(L_i)$ consists of finite order elements (see
e.g. \cite[Lem. 4.7]{BCGP}), so that $T$ is a finite index subgroup
in a product of orbifold surface groups.
This concludes the proof of Proposition \ref{prop:extension}.


\begin{thebibliography}{BCGP09}

\bibitem[Arm65]{armstrong1}
M.~A. Armstrong.
\newblock {\em On the fundamental group of an orbit space},
\newblock {Proc. Cambridge Philos. Soc.} {\bf 61} (1965) 639--646.

\bibitem[Arm68]{armstrong} M. A. Armstrong. \textit{The fundamental
    group of the orbit space of a discontinuous group},
  Proc. Cambridge Philos. Soc. {\bf 64} (1968) 299--301.  

\bibitem[BC04]{newsurfaces}
I. Bauer and F. Catanese,
\newblock \textit{Some new surfaces with {$p_g=q=0$}},
\newblock in {\em The {F}ano {C}onference} (2004), 123--142,
Univ. Torino, Turin.

\bibitem[BCG08]{BCG} I. Bauer, F. Catanese,
F. Grunewald. \textit{The classification of surfaces with {$p_g=q=0$}
isogenous to a product of curves}, Pure Appl. Math. Q. {\bf 4}
(2008), 547--586. 

\bibitem[BCGP09]{BCGP} I. Bauer, F. Catanese, F. Grunewald,
  R. Pignatelli. \textit{Quotients of products of curves, new surfaces
    with $p_g=0$ and their fundamental groups}, arXiv:0809.3420.  

\bibitem[BCP06]{survey}
I. Bauer, F. Catanese, and R. Pignatelli,
\newblock \textit{Complex surfaces of general type: some recent
  progress}, 
\newblock in {\em Global aspects of complex geometry} (2006), 
1--58, Springer, Berlin.

\bibitem[Bea95]{beardon}
A. F. Beardon.
\newblock {\em The geometry of discrete groups}, volume~{\bf 91} of {
  Graduate Texts in Mathematics},
\newblock Springer-Verlag, New York, 1995.
\newblock Corrected reprint of the 1983 original.

\bibitem[Cat00]{fabrizio} F. Catanese.  \textit{Fibred surfaces, varieties
    isogenous to a product and related moduli spaces},
  Amer. J. Math. {\bf 122} (2000), no. 1, 1--44. 

\bibitem[Cat03]{fabrizio2}
F. Catanese.
\newblock {\em Moduli spaces of surfaces and real structures},
\newblock {Ann. of Math.} {\bf 158} (2003) 577--592.


\bibitem[GJZ08]{GJ-ZZ08} F. Grunewald, A. Jaikin-Zapirain,
  P.A. Zalesskii. \textit{Cohomological goodness and the profinite
    completion of {B}ianchi groups}, Duke Math. J. {\bf 144} (2008),
  no. 1, 53--72.

\bibitem[K93]{kollar} J. Koll\'ar. \textit{Shafarevich maps and
    plurigenera of algebraic varieties}, Invent. Math. {\bf 113}
  (1993), no. 1, 177--215.

\bibitem[Ser94]{serre}
J.-P. Serre.
\newblock {\em Cohomologie galoisienne}, volume~{\bf 5} of {Lecture
  Notes in Mathematics},
\newblock Springer-Verlag, Berlin, fifth edition, 1994.



\end{thebibliography}

{\small
\vskip .4cm \noindent
\textsc{%
Universit\'e Paul Sabatier,
Institut de Math\'ematiques de Toulouse,
118 route de Narbonne,
F-31062 Toulouse Cedex 9, France} \\
\texttt{thomas.dedieu@math.univ-toulouse.fr}

\vskip .2cm \noindent
\textsc{%
Lehrstuhl Mathematik VIII, Mathematisches Institut,
Universit\"at Bayreuth, D-95440 Bayreuth, Germany}\\
\texttt{fabio.perroni@uni-bayreuth.de}
}

\end{document}